\documentclass[twocolumn]{autart}    

\usepackage{graphicx}          
\usepackage{amsmath}
\usepackage{amssymb}
\usepackage{color}
\usepackage[mathscr]{euscript}
\renewcommand{\l}{\mathbb {I}}
\newcommand{\rank}{\operatorname{rank}}

\begin{document}

\begin{frontmatter}

\title{Closed-form expressions for the pure time delay in terms of the input and output Laguerre spectra } 

\thanks[footnoteinfo]{This paper was not presented at any IFAC 
meeting.}

\author[Uppsala]{Alexander Medvedev}\ead{alexander.medvedev@it.uu.se}  

\address[Uppsala]{Information Technology, Uppsala University, Uppsala, SWEDEN}  


\begin{abstract}                          
The pure time delay operator is considered in continuous and discrete time under the assumption of the input signal being integrable (summable) with square. By making use of a discrete convolution operator with polynomial Markov parameters, a common framework for handling the continuous and discrete case is set.
Closed-form expressions for the delay value are derived in terms of the Laguerre spectra of the output and input signals. The expressions hold for any feasible value of the Laguerre parameter and can be utilized for e.g. building time-delay estimators that allow for non-persistent  input. A simulation example is provided to illustrate the principle of Laguerre-domain time delay modeling and analysis.
\end{abstract}

\end{frontmatter}

\section{Introduction}
Describing  properties of signals and systems  as functions of time is a natural option since time domain is the most common and easily interpretable observation framework. Yet, it does not necessarily yield the best setup for system analysis, estimation, and design when the involved signals are known to possess a certain property.  For instance,  continuous periodic signals are efficiently represented by Fourier series which description has led to the development of  frequency-domain methods \cite{Z53}. Besides introducing parsimonious signal representations, the use of frequency domain has  simplified the mathematical tools for analysis of dynamical system by substituting   time-domain convolution operators with multiplication of  the corresponding Fourier transforms. 

Ubiquitous in systems theory signals integrable with square  have suggested the use of Laguerre functions \cite{C82} in their modeling. In \cite{NY82}, the role of  Laguerre functions as a versatile tool for studying properties of discrete linear time-invariant systems was convincingly demonstrated. In fact, the Laguerre shift operator applied to generate a Laguerre (functional) basis \cite{MP99} is, in a well-defined sense, equivalent to the forward shift ($q$) operator and can be generally employed for describing signals and systems in both continuous and discrete time. As argued in \cite{MBU20}, the Laguerre shift is one of the bilinear operators, along with $\delta$-operator, $\gamma$-operator, Tustun's operator, proposed in systems theory, primarily to gain beneficial numerical properties of the manipulated objects and unify continuous and discrete frameworks.

A promising but underdeveloped application area of Laguerre functions is time-delay systems. In contrast with the solid results available regarding  approximation  of time-delay systems by means of orthogonal functional bases, see e.g. \cite{MP99}, \cite{MP99a}, publications on time-delay systems in Laguerre domain, i.e. when the involved signals are represented by their Laguerre spectra, are scarce.  References to  relevant papers along this avenue  of research that originates from \cite{FM99} are provided in the next section.

Tme-of-flight estimation of signals with finite energy (pulses) that constitutes the core of radar, sonar, ultrasound, and lidar technology \cite{M92} is essentially the estimation of the delay between an emitted and reflected pulse. Even though delay estimation methods based on Laguerre functions were numerically benchmarked  against conventional techniques in \cite{BL03}, lack of system theoretical grounds still hinders their practical applications. 




This paper  focuses on the mathematical description of the time-delay operator in Laguerre domain.
The main contributions of this work are as stated below.
\begin{itemize}
\item A common framework for describing pure continuous and discrete time-delay operators in  Laguerre domain in the form of a convolution operator with polynomial Markov parameters is introduced. It enables  analysis, design, and estimation approaches uniformly applicable  to  both continuous and discrete  delay models.
\item Based on the properties of the proposed convolution description, closed-form expressions for the time-delay value as a function of the input and output Laguerre spectra   delay block are derived.
\end{itemize}

Notice that closed-form expressions for the delay are not intended as estimators since they do not take into account model uncertainty and disturbance properties.  Yet, as \cite{AM22} demonstrates, discrete  delay estimation    in the face of heavily correlated and non-stationary additive measurement disturbances can be built  in Laguerre domain by exploiting the analytical results. 

The paper is organized as follows. First necessary notions of Laguerre domain are introduced. Then, the Laguerre domain representations of the continuous and discrete pure delay operators are revisited and cast in a common framework. Making use of the latter, the delay value is analytically related to the input and output Laguerre spectra of the delay block. The presented concepts are illustrated by simulation.

\section{Laguerre domain}
\paragraph*{Continuous time:}The Laplace transform of the $k$-th  continuous Laguerre function  is given by
\[
\ell_k(s)=\frac{\sqrt{2p_c}}{s+p_c}T_c^k(s), \quad T_c(s)\triangleq\frac{s-p_c}{s+p_c},
\]
for $k \in \mathbb{N}$, where $p_c>0$ represents the Laguerre parameter, and $T_c$ is the continuous Laguerre shift operator. 

Let $\mathbb{H}_c^2$  be the Hardy space of  functions analytic in the open left half-plane.
The set $\{\ell_k \}_{k \in \mathbb{N}}$ is an orthonormal complete basis in $\mathbb{H}_c^2$  with respect to the inner product 
\begin{align}\label{eq:inner}
\langle W, V \rangle_c \triangleq \frac{1}{2\pi i}\int_{-\infty}^{\infty}W(s)V(-s)~{\rm d}s.
\end{align}

Any function $W \in \mathbb{H}^2_c$  can be represented as a series
$$ W(s)=\sum_{k=0}^\infty w_k \ell_k(s),  \quad  w_j=\langle W, \ell_j \rangle_c ,$$
and the set $\{w_j\}_{j\in\mathbb{N}}$ is then  referred to as the {\em continuous Laguerre spectrum} of $W$.

\paragraph*{Discrete time:}
The  discrete Laguerre functions  are specified in  $\mathscr Z$-domain  by
\begin{equation}\label{eq:Laguerre}
    L_j(z)=\frac{\sqrt{1-p_d}}{z-\sqrt{p_d}}T_d^j(z), \quad T_d(z)\triangleq \frac{1-\sqrt{p_d}z}{z-\sqrt{p_d}}, 
\end{equation}
for all $j \in \mathbb{N}$, where the constant $0<p_d<1$ is the discrete Laguerre parameter, and $T_d(z)$ is the discrete shift operator.

Let $\mathbb{H}_d^2$ be the Hardy space of analytic functions on the complement of the unit disc that are square-integrable on the unit circle and equipped with  the inner product 
\begin{align}\label{eq:inner_disc}
\langle W, V \rangle_d = \frac{1}{2\pi i}\oint_{D} \ W(z)\overline{V(z)}~\frac{{\mathrm d}z}{z},
\end{align}
where ${\overline{V(z)}=V(z^{-1})}$ and $D$ is  the unit circle. Then,  $\{L_k(z) \}_{k \in \mathbb{N}}$  is an orthonormal complete basis in  $\mathbb{H}^2_d$.

 Any function $W \in \mathbb{H}^2_d$  can be represented as a series
$$ W(z)=\sum_{k=0}^\infty w_k L_k(z),  \quad  w_j=\langle W, L_j \rangle_d ,$$
and the set $\{w_j\}_{j\in\mathbb{N}}$ is referred to as the {\em discrete Laguerre spectrum} of $W$.

A system is said to be represented in Laguerre domain when the involved in the system model signals are represented by their Laguerre spectra.

Further on, both the continuous and discrete case are presented in a uniform manner and no notational difference is made when the framework is clearly stated.

\section{Time delay  in Laguerre domain}\label{sec:time_delay}

\paragraph*{Continuous time:}
The well-known associated Laguerre polynomials  (see e.g.  \cite{S39})   are explicitly given by
\begin{equation}\label{eq:assoc_Laguerre}
\mathrm{L}_m(\xi; \alpha)=\sum_{n=0}^m \frac{1}{n!} \binom{m+\alpha}{m-n} (-\xi)^n, \ \forall m\in \mathbb{N}, \quad \xi\in\mathbb{R}.
\end{equation}
In what follows, only the polynomials with a particular value of $\alpha$ are utilized and the shorthand notation $\mathrm{L}_m(\xi)\triangleq\mathrm{L}_m(-\xi; \alpha)|_{\alpha=-1}$ is introduced. 

Consider the signal $u(t) \in \mathbb{L}_2\lbrack 0, \infty)$ given by its Laguerre spectrum $\{ u_j \}_{j\in\mathbb{N}_0}$.  Being passed through to a delay block 
\begin{equation}\label{eq:cdelay}
y(t)=u(t-\tau), \quad \tau\ge 0, \quad t\in \lbrack 0, \infty),
\end{equation}
 the input $u(t)$ results in the output  $y(t) \in \mathbb{L}_2\lbrack 0, \infty)$ with the spectrum $\{ y_j \}_{j\in\mathbb{N}_0}$. Then, according to \cite{Hidayat12}, the following relationship holds  between the spectra
\begin{align}\label{eq:poly_pure_delay}
y_j = \sum_{k=0}^{j-1} h_{j-k}(\kappa) u_{k}+ h_0(\kappa) u_j , \ \forall j \in \mathbb{N}_0,  
\end{align}
where $h_k(\kappa)=\mathrm{e}^{-\frac{\kappa}{2}}\mathrm{L}_k(\kappa)$, and  $ \kappa=2p\tau$.
Notice that $\mathrm{L}_0(\kappa)=1$ and, therefore,  $h_0(\kappa)=\mathrm{e}^{-\frac{\kappa}{2}}$.


\paragraph*{Discrete time:}
Introduce the polynomials
\[
L_m^{(\tau)}(\xi) = (-\xi)^{m-\tau}\sum_{n=0}^{\tau-1}\binom{m+n}{n}\binom{m-1}{\tau-n-1}(-\xi^2)^n,
\]
where it is agreed that $\binom{n}{k}=0$ for $k>n$ by definition.
For the discrete delay operator in time domain, it holds
\begin{equation}\label{eq:ddelay}
y(t)= u(t-\tau), \quad t, \tau\in\mathbb{N}_0.
\end{equation}
Assuming  $u(t) \in \l^2\lbrack 0,\infty)$ and with the same notation for the input and output Laguerre spectra, it is shown in \cite{MBU20} that 
\begin{align} \label{eq:delay-coeffs}
    y_j = \sum_{k=0}^{j-1}
    h_{j-k}(\xi)u_k + h_0(\xi) u_j, \ \forall j \in \mathbb{N}_0,  
    \end{align}
where
$h_0(\xi)=\xi^\tau$, $h_j=(1-p) L_{j}^{(\tau)}(\xi), j\ge 1$, and $\xi=\sqrt{p}$.

A readily observed difference between the continuous and discrete time cases is  that the argument  of the polynomials $\mathrm{L}_k(\kappa)$ carries information about $\tau$ whereas the argument of $L_m^{(\tau)}(\xi)$ is solely defined by the Laguerre parameter $p$. In the continuous-time case of \eqref{eq:poly_pure_delay}, the delay value appears only in a product with the Laguerre parameter thus highlighting the role of the latter as a time scale degree of freedom. On the contrary, in the discrete case of \eqref{eq:poly_pure_delay}, $\tau$ influences the order of $L_m^{(\tau)}(\xi)$ as the time scale is fixed by the discrete time variable $t$.

The common convolution form of the continuous and discrete delay operators, i.e. \eqref{eq:poly_pure_delay}  and \eqref{eq:delay-coeffs}, implies that these descriptions exhibit some kind of ``causality". Indeed, the output coefficient $y_k$ depends only on the input coefficients $u_k, u_{k-1}, \dots, u_0$.  Naturally, this is not a temporal casualty  since each Laguerre coefficient is evaluated on the time interval $t\in \lbrack 0, \infty)$.  Assuming that the input signal is formed so that $u_i=0, i=0,\ldots,m-1$, for some $m\in \mathbb{N}_+$, the first $m$ coefficients of the output, i.e. $y_i, i=0,\ldots,m-1$, are independent of the input and constitute instead the first $m$ coefficients of the Laguerre spectrum  of an additive disturbance. Thus, the signal shape of the realization can be reconstructed (see Section~\ref{sec:ex}) and utilized, e.g.  for noise reduction.

Despite the  identical  form of Laguerre-domain representations of the continuous and discrete delay in \eqref{eq:cdelay} and \eqref{eq:ddelay},  the continuous operator is an infinite-dimensional system and the discrete operator admits a minimal state-space realization of order $\tau$. The difference  roots in the properties of the polynomials $\mathrm{L}_m(\cdot)$ and $L_m^{(\tau)}(\cdot)$. 

Recalling that the output of a discrete linear time-varying system under zero initial conditions is given by the convolution of the input sequence with the system's Markov parameters,  $h_0, h_1, \dots$ are further referred to as the Laguerre-domain Markov parameters of delay operators \eqref{eq:cdelay}, \eqref{eq:ddelay}.

Consider now a Hankel matrix composed of the Markov parameters
\[
\mathcal{H}_{n}=\begin{bmatrix}
h_0 &h_1  &\ldots &h_{n-1} \\
h_1 &h_2   &\ldots &h_n\\
\vdots &\vdots   &\ddots &\vdots\\
h_{n-1} &h_n   &\ldots  &h_{2n-2}\\
\end{bmatrix}.
\]
\begin{prop}\label{th:Markov}
For the Markov parameters in \eqref{eq:poly_pure_delay}
\[ \rank \mathcal{H}_{n}=n, \quad \forall n.
\]
For the Markov paraemters in \eqref{eq:delay-coeffs}
\[ \rank \mathcal{H}_{n}= \begin{cases}  n & \text{ if $n<\tau$},\\
\tau & \text{if $\tau \le n$.}
\end{cases}
\]
\end{prop}
\begin{pf} The Ho-Kalman algorithm \cite{HK66} relates the rank of $\mathcal{H}_{n}$ to the order of the minimal (continuous or discrete) state-space realization (McMillan degree) that yields the Markov parameters. 
\paragraph*{Continuous time:}  In view of the definition of the Markov parameters in \eqref{eq:poly_pure_delay}, it suffices to show that the polynomials $\mathrm{L}_m$ constitute an orthogonal functional basis  on the positive real axis.

In \cite[ p.~205]{R60}, there is a proof for the orthogonality property of the associated Laguerre polynomials $\mathrm{L}_m(\xi; \alpha)$ (as defined in \eqref{eq:assoc_Laguerre})
assuming $\alpha>-1$.  However, the case of $\alpha=-1$ is  not  covered there, while it is exactly the type of Laguerre polynomials that are encountered  in time-delay systems.
Notably, for  $\alpha<-1$, Laguerre polynomials are no longer orthogonal on the positive real axis but rather satisfy a non-Hermitian orthogonality on certain contours in the complex plane, cf. \cite{KM01}.  Some confusion yet arises with respect to what notion of orthogonality applies in case of $\alpha=-1$, cf. \cite[p.~205 and p.~209]{KM01}.  Therefore, orthogonality of  the polynomials $\mathrm{L}_m$ is proven in Appendix~\ref{app:poly} thus securing the fact that the matrix $\mathcal{H}_{n}$ is always full rank for the continuous time delay case.


\paragraph*{Discrete time:}
As shown in \cite{M22}, a  minimal realization of order $\tau$  satisfying \eqref{eq:delay-coeffs} is 
given by the Laguerre-domain state-space equations
\begin{align}\label{eq:dstate_laguerre}
x_{j+1}&= Fx_j +Gu_j,  \\
y_j&= Hx_j+ {Ju_j}, \nonumber
\end{align}
where
\begin{align*}
F&=\begin{bmatrix} 
-\sqrt{p} &1-p &\dots  &\sqrt{p}^{\tau-2}(1-p) \\
0		& -\sqrt{p} & \dots &  \sqrt{p}^{\tau-3}(1-p) \\
\vdots &\vdots &\vdots &\vdots \\
0 & 0 &\dots &-\sqrt{p}
\end{bmatrix}, \\
G^\intercal &=(1-p) \begin{bmatrix}  \sqrt{p}^{\tau-1}  & \sqrt{p}^{\tau-2}  &\dots &1  \end{bmatrix}, \\
H&= \begin{bmatrix} 1 & \sqrt{p} & \dots  & \sqrt{p}^{\tau-1}  \end{bmatrix},
\end{align*}
and  $J=\sqrt{p}^{\tau}$ $\blacksquare$
\end{pf}

The realization in \eqref{eq:dstate_laguerre} reveals that the delay $\tau$ appears in the discrete-time case both as the system order and system parameter. This is in a sharp contrast with the time-domain description, where $\tau$ is the order of the system and not a parameter.  System order estimation techniques  \cite{B01} are less elaborate compared to those of  parameter estimation. Thus, the representation of the pure time-delay operator in Laguerre enables applying the well-developed technology of parameter estimation \cite{AM22}. The next section demonstrates that the delay value can be obtained both in continuous and discrete time without an intermediate Markov parameter representation   and directly from the  input and output spectra.

\section{Main result}

\begin{prop}\label{th:convolution}
Consider a system given by the  Markov parameters $\{ h_0, h_1, \dots \}$, $h_i \in \mathbb{R}, i=0,1, \dots$, driven by the input sequence $ \{ u_0, u_1, \dots \}$, $u_i \in \mathbb{R}, i=0,1, \dots$,  $u_0\ne 0$, and producing the output sequence $ \{ y_0, y_1, \dots \}$, $y_i \in \mathbb{R}, i=0,1, \dots$ according to 
\begin{align}\label{eq:convolution}
    y_j = \sum_{k=0}^{j-1}
    h_{k}u_{j-k} + h_0 u_j.
    \end{align}
 Then, it applies that
\begin{equation}\label{eq:markov_conv}
h_k= \sum_{j=0}^k g_{k-j}y_j,
\end{equation}
and 
\begin{align*}
g_0&= \frac{1}{u_0} , \\
g_k&= - \frac{1}{u_0} \sum_{j=0}^{k-1} u_{k-j} g_j, \quad k\ge 1. 
\end{align*}
\end{prop}
\begin{pf}
Introduce the following notation
\[
Y_N\!\!=\!\!\begin{bmatrix}
y_0 \\ y_1 \\ \vdots \\ y_{N}
\end{bmatrix}; \  T(U_N)
\!\!=\!\!\begin{bmatrix}
u_0 &0 &\dots &0 \\
u_1 & u_0 &\dots &0\\
\vdots &\vdots & \ddots  &0\\
u_{N} & u_{N-1} &\dots &u_0
\end{bmatrix}; \  H_N\!\!=\!\!
\begin{bmatrix}
h_0 \\  h_1 \\ \vdots \\ h_{N}
\end{bmatrix}.
\]
Then, \eqref{eq:convolution} implies that
\begin{equation}\label{eq:system_markov}
Y_N=T(U_N) H_N.
\end{equation}
The assumption $u_0\ne 0$ secures non-singularity of $T(U_N)$ and the Markov coefficients can be uniquely recovered from the Laguerre spectra of the input and output
\[
H_N=T^{-1}(U_N)Y_N.
\]
The result now follows by a direct application of Lemma~\ref{lm:toeplitz} in Appendix~\ref{app:toeplitz}
$\blacksquare$
\end{pf}

The assumption of $u_0\ne 0$ is not restrictive. The leading zero coefficients can be skipped thus implying that corresponding coefficients of the output are also zero, according to the ``causality" property in Laguerre domain.

Now the main results of the paper can be formulated   stating that 
the value of $\tau$ in \eqref{eq:ddelay}  or \eqref{eq:cdelay} can  be analytically evaluated from the input and the output  spectra.

\begin{prop}\label{th:delay_value}
Consider three  subsequent Laguerre-domain Markov parameters  $h_{m-1},h_m, h_{m+1}$ evaluated for any $m\ge 1$ from the Laguerre spectra of the input and output signals  of either \eqref{eq:cdelay} or \eqref{eq:ddelay}  $ \{ u_0, u_1, \dots \}$, $ \{ y_0, y_1, \dots \}$  according to  \eqref{eq:markov_conv}.
  The following relationships hold then for any admissible value of the Laguerre parameter $p$:
\begin{description}
\item[Continuous time:]
\begin{equation}\label{eq:kappa_est}
\kappa = -\frac{(m+1)h_{m+1} + (m-1)h_{m-1}- 2mh_{m}}{h_{m}},
\end{equation}
where $\kappa=2p\tau$.
\item[Discrete time:]
\begin{equation}\label{eq:tau_est}
\tau
= -\frac{ (m+1)h_{m+1}+ (m-1)h_{m-1}  +m (\xi +\xi^{-1})h_{m}}{(\xi -\xi^{-1})h_{m} } ,
\end{equation}
where $\xi=\sqrt{p}$.
\end{description}
\end{prop}
\begin{pf}
From \eqref{eq:poly_pure_delay}  and \eqref{eq:delay-coeffs}, it is readily observed that the Markov parameters $h_0, h_1, \dots$ are proportional to the corresponding polynomials  $\mathrm{L}_m$ and $L_m^{(\tau)}$. Then it suffices to show that  identities \eqref{eq:kappa_est} and \eqref{eq:tau_est} apply to the polynomials.
\paragraph*{Continuous time}
Recall that the Markov parameters of the continuous delay are defined as
$$h_m(\kappa)=\mathrm{e}^{-\frac{\kappa}{2}}\mathrm{L}_m(\kappa)=\mathrm{e}^{-\frac{\kappa}{2}}\mathrm{L}_m(-\kappa;\alpha)_{\alpha=-1}.$$ 
It is well known, see e.g.  \cite[Theorem~68]{R60}, that  the zeros of $\mathrm{L}_m(x;\alpha),\alpha>-1$ are positive and distinct, which excludes the case at hand. Yet, in view of Proposition~\ref{pr:polynomials}, the polynomials $\mathrm{L}_m(\kappa)$ possess orthogonality on $\lbrack 0, \infty)$. Then, in virtue of \cite[Theorem~55]{R60}, all the zeros of $\mathrm{L}_m(\kappa)$  are distinct and negative. Therefore, it is guaranteed that $  h_m(\kappa)\ne 0, m=0,1,\dots$ for all feasible values of $\kappa$ and the denominator of \eqref{eq:kappa_est} does not turn to zero.

The associated Laguerre polynomials obey the following three-term relationship \cite{R60}
\begin{equation}\label{eq:three-term-cont}
\mathrm{L}_{m+1}(\kappa)= \frac{1}{m+1}(\kappa+2m)\mathrm{L}_m(\kappa) - \frac{m-1}{m+1} \mathrm{L}_{m-1}(\kappa).
\end{equation}
Solving \eqref{eq:three-term-cont} with respect to $\kappa$ and substituting the Markov parameters in it gives \eqref{eq:kappa_est}. Notice that the term corresponding to  $\mathrm{L}_{m-1}(\kappa)$ becomes zero for $m=1$.  Therefore, the value of $\mathrm{L}_{2}(\kappa)$ is completely defined by that of $\mathrm{L}_{1}(\kappa)$.  Consequently, $\mathrm{L}_{0}(\kappa)$ and, thus, $h_0$ are immaterial to the recursion.
\paragraph*{Discrete time} To secure that $h_m\ne 0$ in \eqref{eq:tau_est}, it is sufficient to recall that the discrete delay operator  in Laguerre domain admits  a finite-dimensional minimal realization given by 
\eqref{eq:dstate_laguerre}. Then, none of the Markov parameters of it turn to zero due to the eigenvalue-revealing structure of the matrix $F$.

Similarly to the continuous case, for \eqref{eq:ddelay}, the Markov parameters are proportional to the polynomials $L_m^{(\tau)}$ that are subject to the following recursion \cite{MBU20} 
\begin{equation} \label{eq:three-term}
L_{m+1}^{(\tau)}(\xi)= a_m^{(\tau)}(\xi)  L_m^{(\tau)}(\xi) + b_m L_{m-1}^{(\tau)}(\xi)
\end{equation}
holds with
\[ a_{m}^{(\tau)}(\xi)\triangleq a_{m,1}^{(\tau)}\xi+a _{m,2}^{(\tau)}\xi^{-1}
\]
\begin{equation}\label{eq:a_m}
a_{m,1}^{(\tau)}=-\frac{m+\tau}{m+1}, \quad a_{m,2}^{(\tau)}=-\frac{m-\tau}{m+1}, \quad b_m= -\frac{m-1}{m+1}.
\end{equation}
Exactly as in the continuous case, since $b_1=0$, the polynomial $ L_{0}^{(\tau)}(\xi)$ is not included in  the recursion and $ L_{2}^{(\tau)}(\xi)$ is completely defined by $ L_{1}^{(\tau)}(\xi)$ $\blacksquare$
\end{pf}

The fact that the discrete Markov parameters $h_m$ are always nonzero has consequences for the zeros of the polynomials $L_m^{(\tau)}$. The polynomials are definitely not orthogonal, as Proposition~\ref{th:Markov} implies. However, one can prove that all the roots of the equation $L_m^{(\tau)}(\xi)=0$ are real and distinct. Further, they are symmetric with respect to origin and located outside of the interval $(-1,1)$ where the admissible values of $\pm \sqrt{p}$ lie. These are the properties of zeros that one would otherwise expect from an orthogonal polynomial family.

In view of the recursive relationships appearing in the proof above, it is instructive to recall Favard’s theorem \cite{Ch78}. It  secures the orthogonality of a polynomial sequence $\phi_n, n=0,1,\dots $ with respect to  some positive weight function if they satisfy 
\begin{equation}\label{eq:Favard}
\eta \phi_n(\eta)= a_n\phi_{n+1}(\eta)+c_n\phi_n(\eta)+d_n \phi_{n-1}(\eta)
\end{equation}
for some numbers $a_n,b_n,c_n$ where $a_n\ne0$, $c_n\ne 0$. Then, it readily follows that the polynomials $L_m^{(\tau)}(\xi)$ in \eqref{eq:three-term-cont} form an orthogonal basis. Yet, Favard’s theorem does not say how to find the weight function, which problem is resolved in Proposition~\ref{pr:polynomials}.

On the contrary, the three-term relationship in \eqref{eq:three-term} does not comply with the conditions of Favard’s theorem since $a_{m}^{(\tau)}(\xi)$ is not a first-order polynomial in $\xi$. In fact, the polynomials $L_m^{(\tau)}(\xi) $ do not even build a polynomial family since order of the polynomials is defined by both $m$ and $\tau$. The question of whether or not the  polynomials in the discrete case are orthogonal is negatively answered by the existence of finite-dimensional realization \eqref{eq:dstate_laguerre} reproducing the sequence of the Markov parameters, see Proposition~\ref{th:Markov}.

Notice also that  there exist much simpler relations between the delay value and some of the Markov parameters than those specified by Proposition~\ref{th:delay_value}, see \cite{M22}. For instance, in continuous time, it applies that
\[ \tau=- \frac{\ln h_0}{p},
\]
and the corresponding identity in discrete time is 
\[ \tau= \frac{2\ln h_0}{\ln p}.
\]
\section{Numerical example}\label{sec:ex}
\begin{figure}
    \centering
    \includegraphics[width=0.5\textwidth]{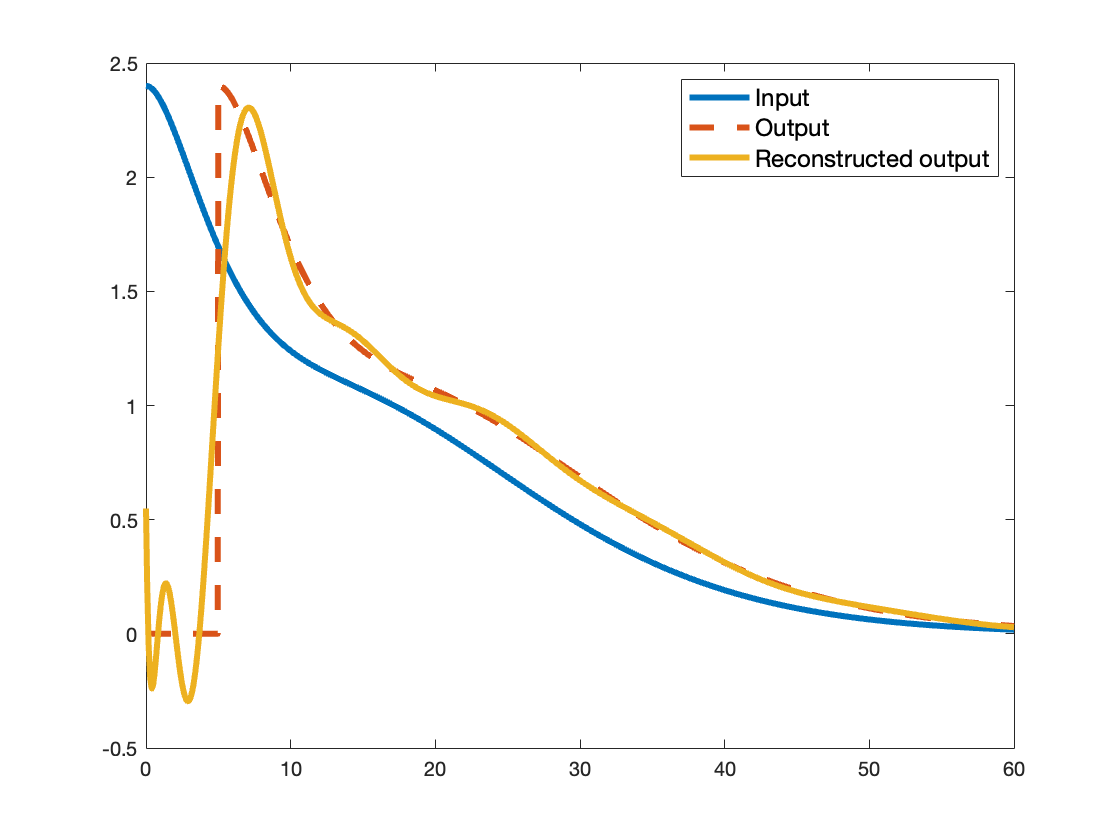}
    \caption{The input signal $u(t)$ of the continuous time-delay operator and its output $y(t)$ (dashed line). An approximation of $y$ by the first 25 Laguerre functions with $p=0.08$ is also shown.}
    \label{fig:io_time}
\end{figure}

\begin{figure}
    \centering
    \includegraphics[width=0.5\textwidth]{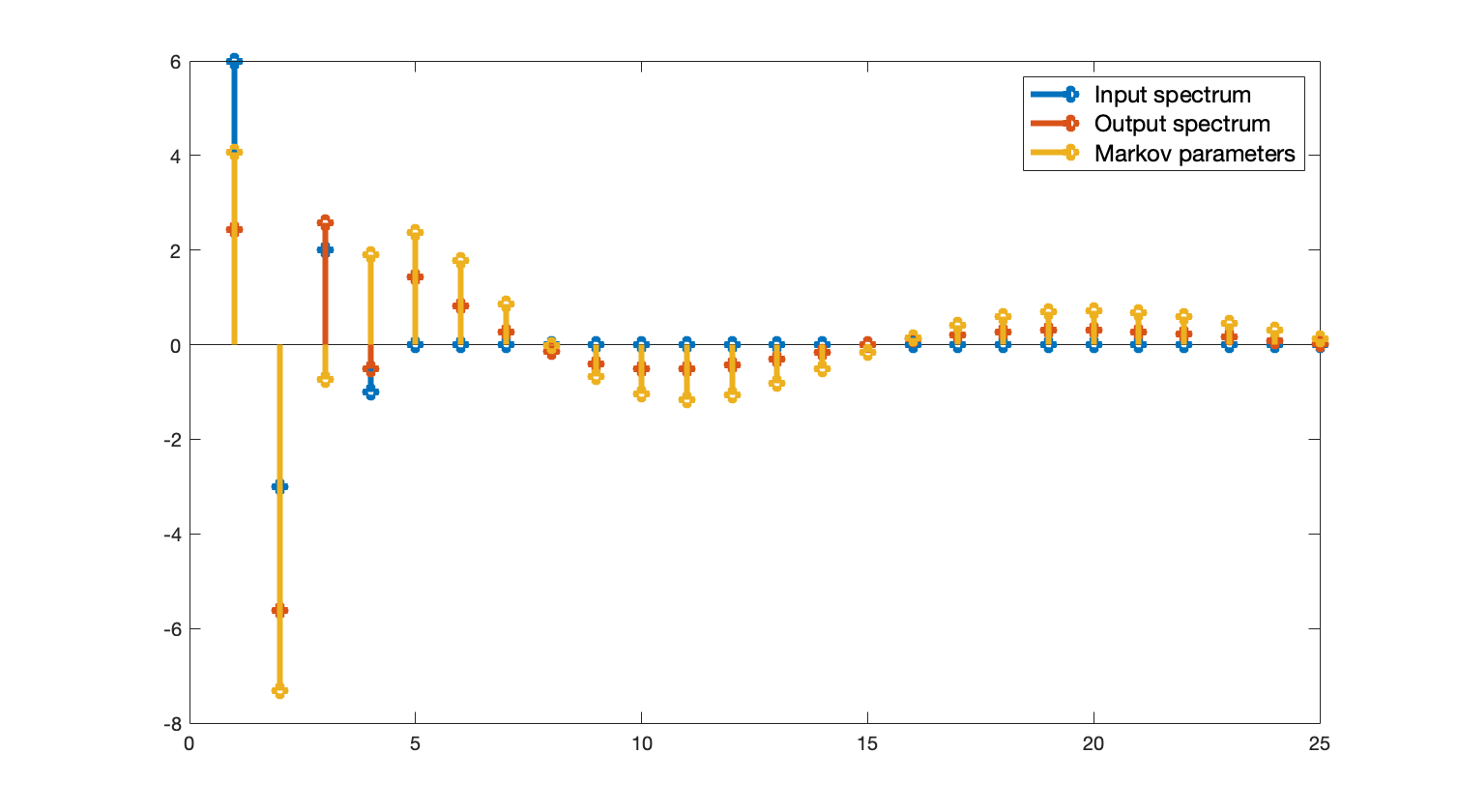}
    \caption{First 25  coefficients of the Laguerre spectra ($p=0.08$) of $u(t)$ and $y(t)$.  Only the four first coefficients of the input spectrum are non-zero. The first 25 Markov coefficients $h_k\times10$ calculated from the input and the output spectra are as well shown.}
    \label{fig:io_spectrum}
\end{figure}

To illustrate the analytical results presented above, consider continuous delay \eqref{eq:cdelay} with $\tau=5$. The input signal is selected to be a linear combination of the first four Laguerre functions
\[
U(s)= 6\ell_0(s)-3\ell_1(s)+2\ell_2(s)-\ell_3(s).
\]
The time evolution of the input as well as the output $y(t)$ are shown in Fig.~\ref{fig:io_time}. The Laguerre spectra of $u(t)$ and $y(t)$ are presented in Fig.~\ref{fig:io_spectrum}. 

Consider now the identity in \eqref{eq:three-term-cont} with $m=2$, i.e.
\begin{equation*}
\kappa = -\frac{3h_{3} + h_{1}- 4h_{2}}{h_{2}},
\end{equation*}
with  $h_1=-0.7318, h_2=  -0.0732, h_3= 0.1903$ (see Fig.~\ref{fig:io_spectrum}) and $p=0.08$. It holds for the numerical values of the example. The same is true for any other feasible $m$.  One should although notice here that the numerical procedures for evaluating Laguerre spectra are not sufficiently developed and the estimates of the coefficients of order over 30 are usually unreliable.

\begin{figure}
    \centering
    \includegraphics[width=0.5\textwidth]{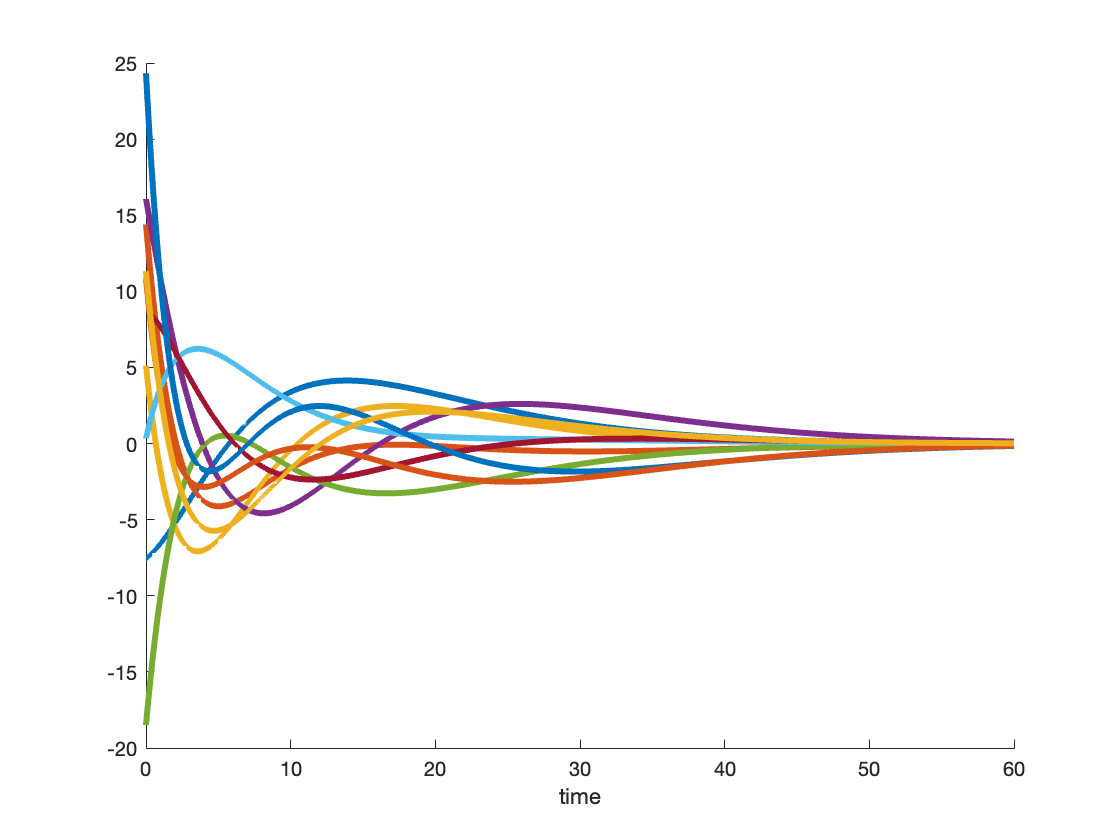}
    \caption{Ten realizations of the disturbance  model representing random linear combination of the first four Laguerre functions, $p=0.18$.  The coefficients of the disturbance model are uniformly distributed random variables in the interval $\lbrack -15, 15\rbrack$.}
    \label{fig:dist}
\end{figure}

\begin{figure}
    \centering
    \includegraphics[width=0.5\textwidth]{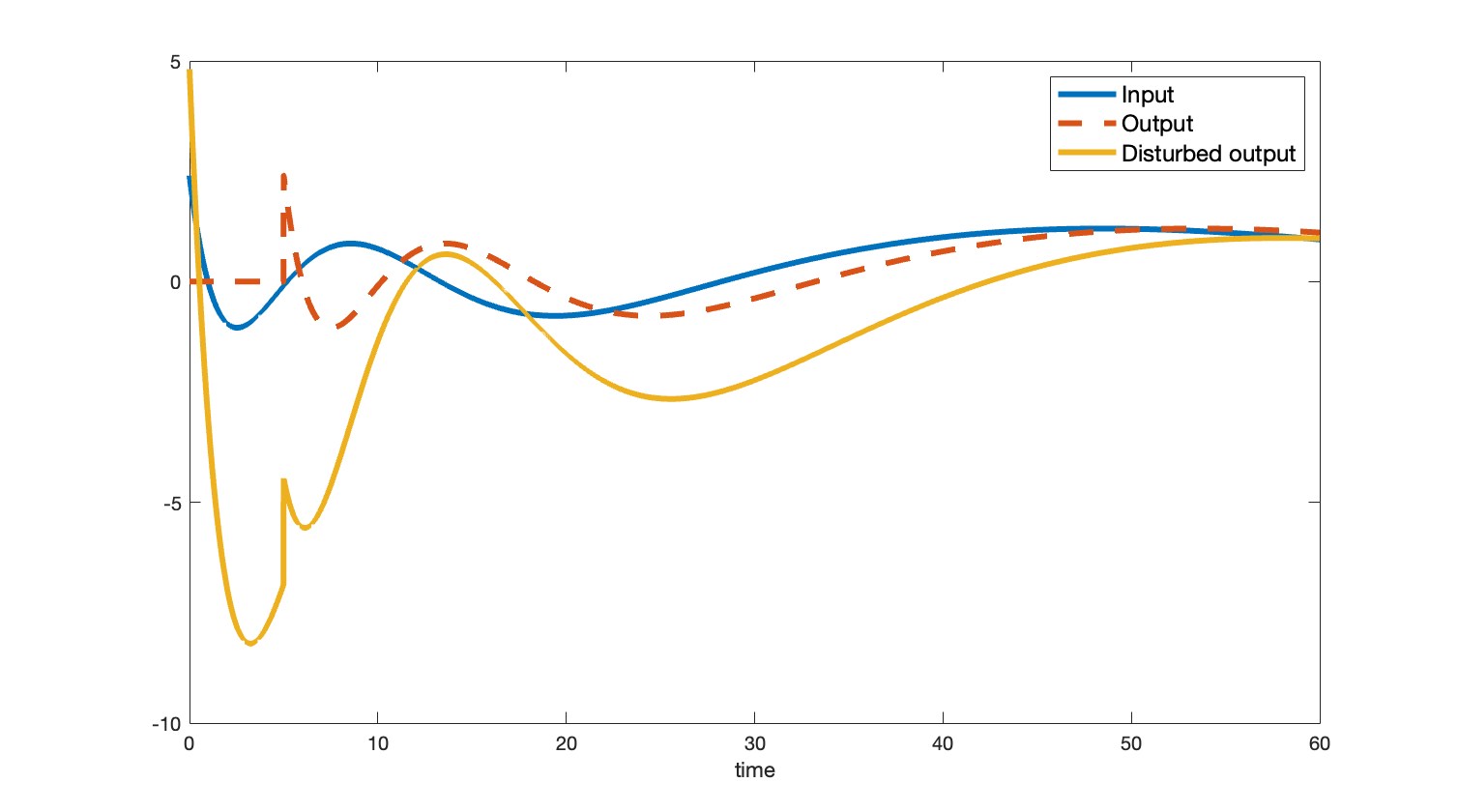}
    \caption{Input signal $u(t)= 6\ell_4(t)-3\ell_5(t)+2\ell_6(t)-\ell_7(t)$, $p=0.18$.  The delay block output $y(t)$, $\tau=5$. The output disturbed with a disturbance realization in Fig.~\ref{fig:dist}.
    The $\mathbb{L}_2$ signal-to-noise ratio is $\| y \|_2 / \| d \|_2 =0.3703$.}
    \label{fig:dist_output}
\end{figure}

\begin{figure}
    \centering
    \includegraphics[width=0.5\textwidth]{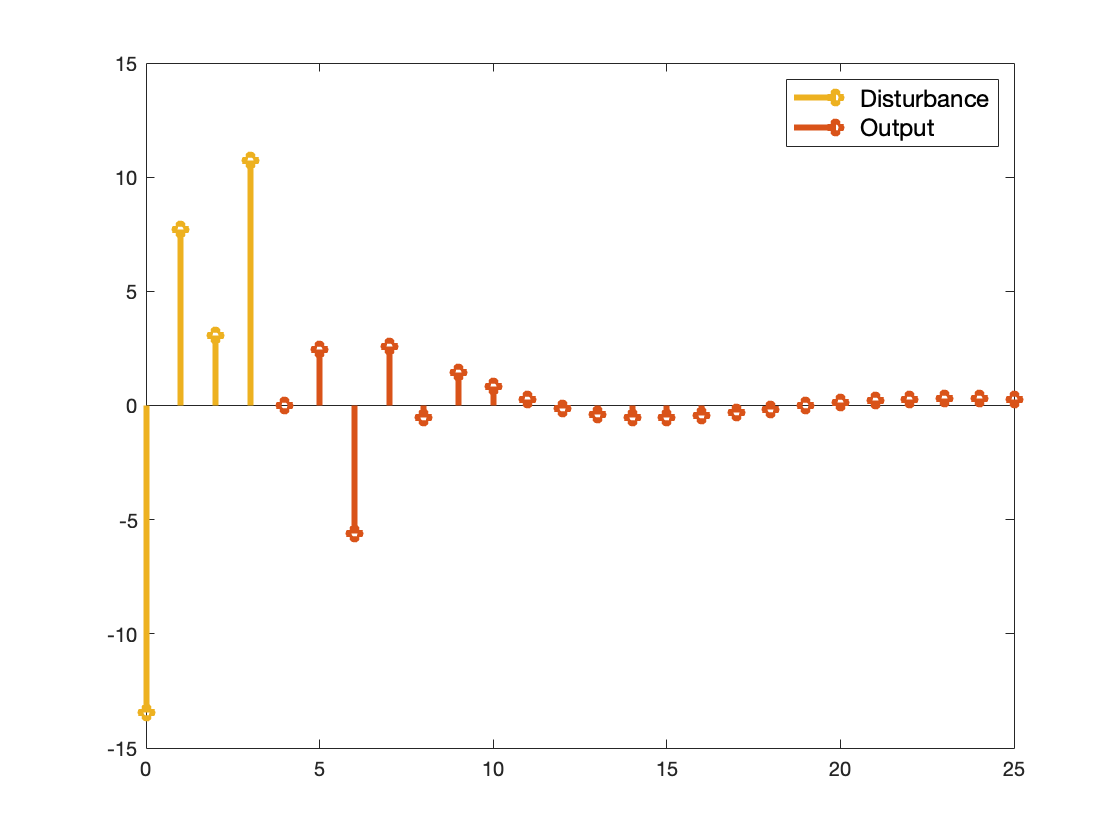}
    \caption{Laguerre spectrum of the disturbed output. The disturbance does not influence the Laguerre coefficients that belong to the delayed  input signal (the delay block output).}
    \label{fig:dist_output_lag}
\end{figure}

As discussed in Section~\ref{sec:time_delay}, both the continuous and discrete descriptions of the time-delay operator exhibit ``casuality" in Laguerre domain. To illustrate the utility of this property, consider a continuous time-delay block with a measurement disturbance
\[ y(t)= u(t-\tau)+d(t),
\]
where $d(t)\in \mathbb{L}_2\lbrack 0, \infty)$ is a random signal. A  Laguerre-domain disturbance model is introduced in \cite{AM22} constituting a linear combination of certain Laguerrre functions with random variables as weights. Fig.~\ref{fig:dist} depicts time-domain realizations of the disturbance model
\[
d(t)= d_1\ell_0(t)+d_2\ell_1(t)+d_3\ell_2(t)+d_4\ell_3(t),
\]
where $d_i, i=1, \dots, 4$ are random numbers uniformly distributed in the interval $\lbrack -15, 15\rbrack$. Now shape the input of the delay block as
\[
u(t)= 6\ell_4(t)-3\ell_5(t)+2\ell_6(t)-\ell_7(t),
\]
so that all the spectral components of $d(t)$ are of lower order compared to those of $u(t)$. Then, in  the Laguerre spectrum of the output signal (see Fig.~\ref{fig:dist_output_lag}), the components that belong to the disturbance do not influence the components of the delayed input which implies that the delay value can be recovered in the same manner than in the disturbance-free case treated above. This is despite the fact that the disturbance completely dominates the output in time domain, Fig.~\ref{fig:dist_output}. Indeed, the disturbance energy is almost three times higher than the energy of the input signal, i.e. $\| d \|_2=19.0958$ and $\| u \|_2=7.0711$. The invariance property is due to the orthogonality of the Laguerre function basis.  A detailed explanation of how the spectral decomposition can be exploited in noise reduction is provided in \cite{AM22}.

\section{Conclusions}

This paper gives closed-form expressions for the delay in terms of the Laguerre spectra of the input and output. Both discrete and continuous time cases are treated. Similarities and differences arising due the infinite-dimensional nature of the continuous delay and a finite-dimensional realization of it in discrete time are highlighted. The derived expressions exhibit the inherent connections between the continuous and discrete delay pure operators by providing a common Laguerre-domain modeling framework exploiting convolution representations with polynomial Markov parameters.

\begin{ack}                               
 This work was partially supported the Swedish Research Council under grant 2019-04451.   
\end{ack}

\bibliographystyle{plain}        
\bibliography{bibliografi}




\appendix
\section{Inverse of lower triangular Toeplitz matrix}\label{app:toeplitz}    
\begin{lem}\label{lm:toeplitz}
If $u_0\ne 0$, the inverse of $T(U_N)$ is given by
\[
T^{-1}(U_N)\triangleq G(U_N)
=\begin{bmatrix}
g_0 &0 &\dots &0 \\
g_1 & g_0 &\dots &0\\
\vdots &\vdots & \ddots  &0\\
g_{N} & g_{N-1} &\dots &g_0
\end{bmatrix},
\]
where 
\begin{align}
g_0&= \frac{1}{u_0} \nonumber \\
g_k&= - \frac{1}{u_0} \sum_{j=0}^{k-1} u_{k-j} g_j, \quad k\ge 1. \nonumber
\end{align}

\end{lem}
\begin{pf}

Introduce the matrix $E_N\in \mathbb{R}^{N\times N}$
$$
E_N=\begin{bmatrix} 0 &0 &\dots &0  &0\\
                    1 & 0 &\dots &0  &0\\
                    \vdots &\vdots &\ddots &\vdots  &0\\
                    0 &0 &\dots &1  &0
\end{bmatrix}.
$$

Then the low-triangular Toeplitz matrix $T(U_N)$ can be written as a matrix-valued polynomial in $E_N$
\[
T(U_N)= \sum_{k=0}^{N-1} u_k E_N^k.
\]
Since $T^{-1}(U_N)$ is also a  low-triangular Toeplitz matrix matrix, it follows
\[
T^{-1}(U_N)= \sum_{k=0}^{N-1} g_k E_N^k.
\]
Observing that $E_N$ is nilpotent
\begin{align}
T(U_N)G(U_N)=G(U_N)T(U_N)&=\sum_{k=0}^{N-1} u_k E_N^k \sum_{k=0}^{N-1} g_k E_N^k \nonumber\\
&= \sum_{k=0}^{N-1} \sum_{i+j=k} u_ig_j E_N^k =I. \label{eq:hg}
\end{align}
The inner sum is taken over all partitions of $k$ into two parts, i.e. the indices $i$ and $j$. For the equality to hold, one has to ensure that
$u_0g_0=1$ and
\[
\sum_{i+j=k} u_ig_j =0, \quad  k\ge 1.
\]
Singling out the term with $g_k$ results in
\[
u_0g_k = - ( u_1g_{k-1} + \dots + u_k g_0    ).
\]
Solving the equation with respect to $g_k$ completes the proof. $\blacksquare$
\end{pf}

\section{Orthogonality of $\mathrm{L}_m$}\label{app:poly}       

\begin{prop}\label{pr:polynomials}
Laguerre polynomials $\mathrm{L}_n(\cdot), n\in \mathbb{N}_0$ are orthogonal on the positive real axis in the sense that
\begin{equation} \label{eq:orth}
\int_0^\infty \frac{ \e^{-x}}{x} \mathrm{L}_n(x)\mathrm{L}_m(x)~\mathrm{d}x= \begin{cases}
0,  & n\ne m\\
\frac{1}{n!}, & n=m.
\end{cases}
\end{equation}
\end{prop}
\begin{pf} 
In \cite{R60}, it is shown that
\begin{eqnarray*}
\lefteqn{(m-n)\int_a^b x^\alpha \e^{-x} \mathrm{L}_n(x; \alpha) \mathrm{L}_m(x; \alpha)~\mathrm{d}x=}\\
& & x^{k+1}\e^{-x} 
\left( \mathrm{L}_m(x;\alpha )\dot{\mathrm{L}}_n(x; \alpha)-\mathrm{L}_n(x; \alpha)\mathrm{\dot{L}}_m(x; \alpha) \right)\biggr\rvert_a^b\biggl.
\end{eqnarray*}
Specialized to $\alpha=-1$, the integration result reads
\begin{equation*}
\varphi_{m,n}(x)=\e^{-x} \left( \mathrm{L}_m(x)\dot{\mathrm{ L}}_n(x)-\mathrm{L}_n(x)\dot{ \mathrm{L}}_m(x) \right)
\end{equation*}
Clearly, due to the exponential factor
\begin{equation} \label{eq:lim_infty}
 \lim_{x\to\infty}\varphi_{m,n}(x) \to 0.
  \end{equation}
Furthermore, the least order term  in $\mathrm{L}_m(x),\ m>0$ is always $x$. Therefore, the least order term  in $\dot{ \mathrm{L}}_m(x),\ m>0$ is always $1$. The least order term in the product $\mathrm{L}_m(x)\dot{ \mathrm{ L}}_n(x),\  n\ne 0, m\ne 0, m\ne n$ is, once again, $x$. Then $\varphi_{m,n}(0)=0$ which fact, together with (\ref{eq:lim_infty}) implies orthogonality of $\mathrm{L}_n(\cdot)$ and $\mathrm{L}_m(\cdot)$ for $m\ne n$, i.e.
\eqref{eq:orth}. It can also be shown that for $m=n$ 
\begin{align*}
\int_0^{\infty} \frac{\e^{-x}}{x}\mathrm{L}_n(x)\mathrm{L}_n(x) ~\mathrm{d}x = \frac{1}{n!}.
\end{align*}
$\blacksquare$
\end{pf}

\end{document}